\documentclass[journal]{IEEEtran}
\usepackage{algorithmic}

\newsavebox{\ieeealgbox}

\let\labelindent\relax   
\usepackage{enumitem}
\setlist[itemize]{leftmargin=*}
\usepackage{graphicx}
\usepackage{array}
\usepackage{textcomp}
\usepackage{cite}
\usepackage{amsmath}
\usepackage{multirow}
\usepackage{bm}
\usepackage{verbatim}
\usepackage{xcolor}

\usepackage{threeparttable}
\usepackage{url}

\newcommand{\splitatcommas}[1]{%
	\begingroup
	\begingroup\lccode`~=`, \lowercase{\endgroup
		\edef~{\mathchar\the\mathcode`, \penalty0 \noexpand\hspace{0pt plus 1em}}%
	}\mathcode`,="8000 #1%
	\endgroup
}

\ifCLASSOPTIONcompsoc
\usepackage[caption=false,font=normalsize,labelfont=sf,textfont=sf]{subfig}
\else
\usepackage[caption=false,font=footnotesize]{subfig}
\fi

\ifCLASSINFOpdf
\else
\fi
\begin{document}
%
\title{The Allocation of a Variable Series Reactor Considering AC Constraints and Contingencies}
%
%
%

\author{Xiaohu~Zhang,~\IEEEmembership{Student Member,~IEEE,}
        Chunlei~Xu, 
        Di~Shi,~\IEEEmembership{Senior Member,~IEEE,} \\
        Zhiwei~Wang,~\IEEEmembership{Member,~IEEE,}
        Qibing~Zhang,
        Guodong~Liu,~\IEEEmembership{Member,~IEEE,} \\
        Kevin~Tomsovic,~\IEEEmembership{Fellow,~IEEE,}
        and Aleksandar~Dimitrovski,~\IEEEmembership{Senior Member,~IEEE}
\thanks{This work was supported in part by ARPAe (Advanced Research Projects
	Agency Energy), in part by the Engineering Research Center Program
	of the National Science Foundation and the Department of Energy under
	NSF Award Number EEC-1041877 and the CURENT Industry Partnership Program, and in part by state grid corporation of China (SGCC).}
\thanks{Xiaohu Zhang, Di Shi and Zhiwei Wang are with GEIRI North America, San Jose, CA, USA, email: \{xiaohu.zhang,di.shi,zhiwei.wang\}@geirina.net. Chunlei Xu and Qibing Zhang are with State Grid Jiangsu Electric Power Company, Naijing, China. Guodong Liu is with the Power and Energy group, Oak Ridge National Laboratory, TN, USA, email: liug@ornl.gov. Kevin Tomsovic is with the Department of Electrical Engineering and Computer Science, the University of Tennessee, Knoxville, TN, USA, email: tomsovic@utk.edu. Aleksandar Dimitrovski is with the Department of Electrical and Computer Engineering, University of Central Florida, FL, USA, email: Aleksandar.Dimitrovski@ucf.edu.}}

%
%

\markboth{This paper is accepted by CSEE JOURNAL OF POWER AND ENERGY SYSTEM, DOI will be provided when it is available.}%
{Shell \MakeLowercase{\textit{et al.}}: Bare Demo of IEEEtran.cls for Journals}
%



\maketitle

\begin{abstract}
The Variable Series Reactors (VSRs) can efficiently control the power flow through the adjustment of the line reactance. When they are appropriately allocated in the power network, the transmission congestion and generation cost can be reduced. This paper proposes a planning model to optimally allocate VSRs considering AC constraints and multi-scenarios including base case and contingencies. The planning model is originally a non-convex large scale mixed integer nonlinear program (MINLP), which is generally intractable. The proposed Benders approach decomposes the MINLP model into a mixed integer linear program (MILP) master problem and a number of nonlinear subproblems. Numerical case studies based on IEEE 118-bus demonstrate the high performance of the proposed approach.
\end{abstract}


%
\IEEEpeerreviewmaketitle


\section*{Nomenclature}
\subsection*{Indices}
\addcontentsline{toc}{section}{Nomenclature}
\begin{IEEEdescription}[\IEEEusemathlabelsep\IEEEsetlabelwidth{$V_1,V_2,V_3$}]
	\item[$i, \ j$] Index of buses.
	\item[$k$] Index of transmission elements.
	\item[$n,\ m$] Index of generators and loads.
	\item[$c$] Index of states; $c=0$ indicates the base case; $c>0$ is a contingency state.
	\item[$t$] Index of load levels.
\end{IEEEdescription}

\subsection*{Variables}
\addcontentsline{toc}{section}{Nomenclature}
\begin{IEEEdescription}[\IEEEusemathlabelsep\IEEEsetlabelwidth{$V_1,V_2,V_3,V_4,V$}]
	\item[$P^g_{nct},Q^g_{nct}$] Active and reactive power generation of generator $n$ for state $c$ under load level $t$.
	\item[$\Delta P^d_{mct}, \Delta Q^d_{mct}$] Active and reactive power load shedding amount of load $m$ for state $c$ under load level $t$.
	\item[$\Delta P^{g,up}_{nct},\Delta P^{g,dn}_{nct}$] Active power generation adjustment up and down of generator $n$ for state $c$ under load level $t$.
	\item[$V_{ict},\theta_{ict}$] Voltage magnitude and angle at bus $i$ for state $c$ under load level $t$.
	\item[$\theta_{kct}$] Voltage angle difference of branch $k$ for state $c$ under load level $t$.
	\item[$x^V_{kct}$] Reactance of a VSR at branch $k$ for state $c$ under load level $t$.
	\item[$\delta_{k}$] Binary variable associated with placing a VSR on branch $k$.
\end{IEEEdescription} 

\subsection*{Funtions}
\addcontentsline{toc}{section}{Nomenclature}
\begin{IEEEdescription}[\IEEEusemathlabelsep\IEEEsetlabelwidth{$V_1,V_2,V_3,V_4,V$}]
	\item[$P_{kct}(\cdot),Q_{nct}(\cdot)$] Active and reactive power flow across branch $k$ for state $c$ under load level $t$..
\end{IEEEdescription} 

\subsection*{Parameters}
\addcontentsline{toc}{section}{Nomenclature}
\begin{IEEEdescription}[\IEEEusemathlabelsep\IEEEsetlabelwidth{$V_1,V_2,V_3,V_4,V$}]
	\item[$r_k,x_{k}$] Resistance and reactance for branch $k$.
	\item[$P^{g,\min}_{nct},Q^{g,\min}_{nct}$] Minimum active and reactive power output of generator $n$ for state $c$ under load level $t$.
	\item[$P^{g,\max}_{nct},Q^{g,\max}_{nct}$] Maximum active and reactive power output of generator $n$ for state $c$ under load level $t$.
	\item[$P_{mct}^d,Q_{mct}^d$]  Active and reactive power consumption of demand $m$ for state $c$ under load level $t$.
	\item[$S_{kct}^{\max}$] Thermal limit of branch  $k$ for state $c$ under load level $t$. 
	\item[$x_{k}^{V,\min},x_{k}^{V,\max}$] Minimum and maximum reactance of the VSR at branch $k$.
	\item[$\theta_{k}^{\max}$] Maximum angle difference across branch $k$.
	\item[$\theta_{i}^{\max},\theta_{i}^{\min}$] Maximum and minimum bus angle at bus $i$.
	\item[$V_{ict}^{\max},V_{ict}^{\min}$] Maximum and minimum bus voltage magnitude at bus $i$ for state $c$ under load level $t$.
	\item[$N_{kct}$] Binary parameter associated with the status of branch $k$ at state $c$ under load level $t$.
	
	\item[$R^{g,up}_n,R^{g,dn}_n$] Ramp up and down limit for generator $n$.
	\item[$a_n^g$] Cost coefficient for generator $n$.
	\item[$a_n^{g,up},a_n^{g,dn}$] Cost coefficient for generator $n$ to increase and decrease active power.
	\item[$a_{LS}$] Cost coefficient for the load shedding.
	\item[$A_h$] Annual operating hours: 8760 h.
	\item[$A^I_k$] \textcolor{black}{Investment cost for a VSR at branch $k$.}
	\item[$\tilde{A}^I_k$] \textcolor{black}{Annualized investment cost for a VSR at branch $k$.}
\end{IEEEdescription} 

\subsection*{Sets}
\addcontentsline{toc}{section}{Nomenclature}
\begin{IEEEdescription}[\IEEEusemathlabelsep\IEEEsetlabelwidth{$V_1,V_2,V_3$}]
	\item[$\mathcal{D},\mathcal{G}$] Set of loads and generators.
	\item[$\mathcal{D}_i,\mathcal{G}_{i}$] Set of loads and generators located at bus $i$.
	\item[$\Omega_{L}$] Set of transmission lines.
	\item[$\Omega_{L}^i$] Set of transmission lines connected to bus $i$.
	\item[$\Omega_{T}$] Set of load levels.
	\item[$\Omega_{c}$] Set of contingency operating states.
	\item[$\Omega_{0}$] Set of base operating states.
	\item[$\Omega_{V}$] Set of candidate transmission lines to install VSR.
	\item[$\mathcal{B}$] Set of buses. 
	\item[$\mathcal{G}_{re}$] Set of generators allowed to rescheduling.
\end{IEEEdescription}
Other symbols are defined as required in the text.

\section{Introduction}
\IEEEPARstart{I}n recent years, with the advent of power market deregulation, the massive integration of the renewable energy and the increasing demand of electricity consumption, the aging power grid has become congested and is under stress, which results in higher energy cost \cite{mybibb:ps_economics,mybibb:national_grid}. Building new transmission lines on some critical corridors is one method to reduce the congestion and increase the system reliability, but the political and environmental constraints make this option unattractive \cite{mybibb:tep_pst_myself}. Therefore, there are growing interests for utilities to actively control the existing power network to manage power flows while at the same time improve security margins and increase the system transfer capability    
\cite{mybibb:optfacts,mybibb:GAFACTS}.

Flexible AC transmission systems (FACTS) are one technology for controlling power flow and enhancing the utilization of existing transmission network \cite{mybibb:FACTS1,mybibb:FACTS2,mybibb:SSSC,mybibb:flow_control_FACTS1}. Specific types of series FACTS devices, which are named as Variable Series Reactor (VSR), have the ability to efficiently regulate the power flow through the adjustment of the transmission line reactance. Typical examples of VSR are Thyristor Controlled Series Compensator (TCSC), Distributed Series Reactor (DSR) and smart wire \cite{mybibb:ETH_TCSC1,mybibb:dsr_pd,mybibb:smart_wire}. According to the Green Electricity Network Integration program (GENI) \cite{mybibb:geni}, more FACTS-like devices \cite{mybibb:Aleks} with far cheaper price will be commercially available soon for the transmission network across US. Hence, the development of efficient algorithms which are capable of finding the optimal locations of FACTS device is of great importance.         

In the technical literature, various strategies have been proposed to allocate and utilize FACTS devices. Due to the nonlinear and non-convex characteristics of the optimal placement problem, various evolutionary computation techniques, such as, genetic algorithm (GA) \cite{mybibb:GAFACTS,mybibb:GA_FACT_market}, differential evolution \cite{mybibb:differentialevoluation}, particle swarm optimization (PSO) \cite{mybibb:PSO_FACTS}, have been proposed to find the optimal locations of TCSC. These techniques have the advantage of straightforward implementation but they do not provide any indicator about the quality of the solutions. Reference \cite{mybibb:TCSC_ali} introduces an index called the single contingency sensitivity (SCS), which provides an indicator for the effectiveness of a given branch in relieving the congestions under all considered contingencies. After the locations of TCSC are selected based on the ranking of SCS, an optimization problem is formulated to obtain the settings of TCSCs for each contingency. To enhance the transfer capability of the network, the authors in \cite{mybibb:TCSCsens} compute the sensitivity of the transfer capability with respect to the line reactance so as to allocate TCSCs. In \cite{mybibb:TCSC_cost_recovery}, sequential optimal power flows are adopted to find the optimal placements of TCSC. The approach is based on a number of optimal power flow (OPF) results from different TCSC locations and settings through a step by step manner. The optimal locations and settings of TCSC are the best optimization results among these OPF results.  

With the rapid development of the branch-and-bound algorithms, the mixed integer program (MIP) has also been employed to solve the power system planning problem \cite{mybibb:tep_es_conje,mybibb:bilevel_wind_tep_csee}.  Reference \cite{mybibb:TCSCLFB} applies the line flow equations \cite{mybibb:LFB} to allocate TCSC. The problem is formulated as a mixed integer linear program (MILP) or mixed integer quadratic program (MIQP). The non-convex bilinear terms in the constraints are eliminated by replacing one variable with its respective limit. However, the limits of these variables, such as, active power flow and active power loss on the transmission lines, cannot be determined a priori, which restricts the utility of the approach. In \cite{mybibb:asu_FACTS1,mybibb:Tao_Ding}, to evaluate the benefits of the VSR on the economic dispatch (ED) problem, the nonlinear term of the product between the variable susceptance and voltage angle is linearized with the big-M method. The original nonlinear program model is transformed into an MILP model and solved by commercial MIP solvers. In \cite{mybibb:tcsc_ac}, the allocation of TCSCs considering load variability is investigated by using Benders Decomposition. The complete model is decomposed into an MILP model that serves as master problem and a number of nonlinear programs (NLP) that serve as subproblems. However, the contingency constraints are not considered.

It has been shown in \cite{mybibb:opf_redispatch_facts} that the series FACTS have the capability of reducing generator rescheduling and load shedding amount following contingencies. Hence, if the contingency constraints are included in the planning model, a more useful investment strategy can be achieved by the system designers. The researchers in \cite{mybibb:PSO_SQP_FACTS} propose a two level hybrid PSO/SQP algorithm to address this problem. The upper level problem is to leverage the standard PSO to determine the locations and capacities of the FACTS devices and the lower level is to decide the settings of the devices for normal state and contingencies by sequential quadratic programming. However, the computational time is high even for small scale system considering limited number of contingencies.   

This paper addresses the optimal placement of VSR in a transmission network considering AC constraints and a series of $N-1$ contingencies. A single target year is considered and three load patterns which represent the peak, normal and low load level are selected to accommodate the yearly load profile. In addition, the $N-1$ contingencies have the probability to occur in any of the load levels. The planning model is a large scale mixed integer nonlinear program (MINLP) model and is quite difficult to be solved by the existing commercial solvers. Therefore, we use Generalized Benders Decomposition to separate the planning model into master problem and a number of subproblems including both the base and contingency operating states. The contribution of this paper are twofold. First, we propose a planning model to optimally allocate the VSR in the transmission network while retaining AC constraints and considering a series of contingencies. Second, we implement a Benders algorithm to solve the proposed model which significantly relieves the computational burden and makes it a potential approach for practical large scale system. 

The remaining sections are organized as follows. In Section \ref{static_model}, the static model of VSR is presented. Section \ref{problem_formulation} illustrates the detailed formulation of the planning model. The solution procedure based on Generalized Benders Decomposition is demonstrated in Section \ref{solution_approach}. In Section \ref{case_study}, the IEEE 118-bus system is selected for case studies. Finally, some conclusions are given in Section \ref{conclusions}.

\section{Static Model of VSR}
\label{static_model}
\begin{figure}[!htb]
	\centering
	\includegraphics[width=0.35\textwidth]{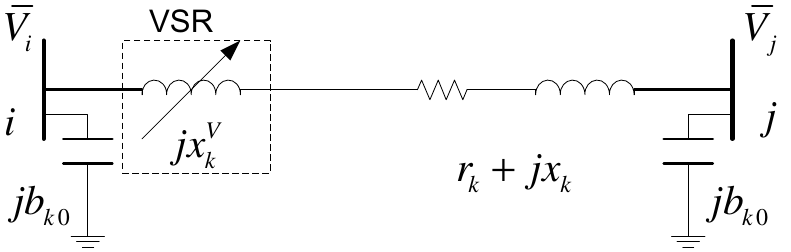}
	\caption{Static representation of VSR.}
	\label{CVSR_steady}
\end{figure}

The static model of VSR can be represented by a variable reactance with the parasitic resistance ignored as given in Fig. \ref{CVSR_steady}. The inserted reactance effectively changes the overall impedance of the branch.

For the candidate lines to install VSR, the total conductance and susceptance of the line should be modified by the following two equations:
\begin{align}
	g'_{k}&=\frac{r_{k}}{r_{k}^2+(x_{k}+\delta_kx_{k}^{V})^2}, \ \ \ k\in \Omega_{V}  \\
	b'_{k}&=-\frac{(x_{k}+\delta_kx^{V}_{k})}{r_{k}^2+(x_{k}+\delta_kx^{V}_{k})^2}, \ \ \ k\in \Omega_{V}
\end{align}

$\delta_k$ is a binary variable which flags the installation of  VSR on line $k$. Hence, the power flow on line $k$ is not only the function of bus voltage magnitude $\bm{V}$ and bus angle $\bm{\theta}$ but also the function of $\bm{x_k^{V}}$ and $\bm{\delta_k}$.

\section{Problem Formulation}
\label{problem_formulation}
The complete optimization model renders a non-convex large scale MINLP. A detailed description of the optimization model will be given in this section.
\subsection{Objective Function}
The objective function employed in this paper is to minimize the total planning cost for the target year. The cost includes three components: 1) \textcolor{black}{annualized investment cost for VSR;} 2) operation cost for normal states; 3) operation cost for contingency states. The objective is then formulated as:
\begin{equation}
	\min_{\Xi_{\text{OM}}} \ \ \sum_{k\in \Omega_V}\tilde{A}_k^I\delta_k+\sum_{t\in \Omega_T}(\pi_{0t}C_{0t}+\sum_{c\in \Omega_c}\pi_{ct}C_{ct}) \label{objective}
\end{equation}  

In (\ref{objective}), $C_{0t}$ is the hourly operating cost under normal states for load level $t$. Assuming a price inelastic load, minimizing the operating cost corresponds to minimizing the generation cost. Then $C_{0t}$ can be expressed as:
\begin{equation}
	C_{0t}=\sum_{n\in \mathcal{G}}a_{n}^gP^g_{n0t}  \label{obj_base}
\end{equation}     
Note that the linear cost coefficient $a_n^g$ is adopted in the objective function. However, if the quadratic cost curve is required for the generators, the piecewise linearization can be used to linearize the cost curve and easily embedded into the model \cite{mybibb:guodong_xiaohu1}.

The hourly operating cost under contingency state $c$ and load level $t$ is denoted as $C_{ct}$, which can be expressed as:
\begin{align}
	C_{ct}&=\sum_{n\in \mathcal{G}}a_{n}^gP^g_{nct}+\sum_{m\in \mathcal{D}}a_{LS}\Delta P^d_{mct} \nonumber \\ 
	&+\sum_{n\in \mathcal{G}_{re}}(a_n^{g,up}\Delta P^{g,up}_{nct}+a_n^{g,dn}\Delta P^{g,dn}_{nct})                        \label{obj_cont} 
\end{align} 
\textcolor{black}{The operating cost during contingencies is categorized into four terms. The first term in (\ref{obj_cont}) corresponds to the generation cost for each contingency; the second term is the load curtailment cost; the  third and fourth term represent the generator rescheduling cost, which indicates that there is a payment to the involving agents on any changes from the base operating conditions \cite{mybibb:opf_redispatch_facts}.} The duration time associated with each operating state is indicated by $\pi_{ct}$. For a single target year, the number of total operating hours is 8760:
\begin{equation}
	\sum_{t\in \Omega_T}\pi_{0t}+\sum_{t\in \Omega_T}\sum_{c\in \Omega_c}\pi_{ct}=A_h    \label{operating_hour}
\end{equation}   

\subsection{Constraints}
The complete set of constraints are given from (\ref{budget})-(\ref{shedding_relation}):
\begin{align}
& \sum_{k \in \Omega_V} \tilde{A}_k^I\delta_k \le A^{\max} \label{budget} \\
&\sum_{n \in \mathcal{G}_i}P^g_{nct}-\sum_{m \in \mathcal{D}_i}(P^d_{mct}-\Delta P^d_{mct}) \nonumber \\
&\ \ \ \ \ \ \ \ \ \ \ \ \ \ \ \ =\sum_{k\in \Omega_{L}^{i}}  N_{kct}P_{kct}(\bm{\theta,V,x^V,\delta}) \label{p_bal}  \\
&\sum_{n \in \mathcal{G}_i}Q^g_{nct}-\sum_{m \in \mathcal{D}_i}(Q^d_{mct}-\Delta Q^d_{mct})  \nonumber \\
&\ \ \ \ \ \ \ \ \ \ \ \ \ \ \ \ \ =\sum_{k\in \Omega_{L}^{i}}  N_{kct}Q_{kct}(\bm{\theta,V,x^V,\delta}) \label{q_bal}  \\
&\Delta P^d_{m0t}=0, \Delta Q^d_{m0t}=0  \label{lS_normal} \\
&N_{kct}\sqrt{P^2_{kct}(\bm{\theta,V,x^V,\delta})+Q^2_{kct}(\bm{\theta,V,x^V,\delta})} \le S_{kct}^{\max} \label{thermal_limit} 
\end{align}
\begin{align}
& P_{nct}^{g,\min}\le P_{nct}^g \le P_{nct}^{g,\max} \label{Pg_limit1} \\
& Q_{nct}^{g,\min}\le Q_{nct}^g \le Q_{nct}^{g,\max} \label{Qg_limit1} \\
& V_{ict}^{\min} \le V_{ict} \le V_{ict}^{\max}  \label{v_limit}   \\
& \theta_{i}^{\min} \le \theta_{ict} \le \theta_{i}^{\max}  \label{theta_limit}  \\
& x_{k}^{V,\min} \le  x_{kct}^V \le  x_{k}^{V,\max}   \label{vsr_limit} \\ 
& \theta_{ref}=0   \label{ref_angle}   \\
&P_{nct}^g=P_{n0t}^g+\Delta P_{nct}^{g,up}-\Delta P_{nct}^{g,dn}, \ n\in \mathcal{G}_{re}    \label{ramp_gen}  \\
&0 \le \Delta P_{nct}^{g,up} \le R_n^{g,up}, \ n\in \mathcal{G}_{re}   \label{ramp_up}   \\
&0 \le \Delta P_{nct}^{g,dn}  \le R_n^{g,dn},\ n\in \mathcal{G}_{re}  \label{ramp_dn} \\
&P_{nct}^g=P_{n0t}^g, \ n\in \mathcal{G}\backslash\mathcal{G}_{re}  \label{g_fix}  \\
&0 \le \Delta P^d_{mct} \le P^d_{mct} \label{load_shedding}  \\
& \Delta P^d_{mct}Q^d_{mct}=\Delta Q^d_{mct}P^d_{mct}  \label{shedding_relation}
\end{align}
Constraints (\ref{p_bal})-(\ref{ref_angle}) hold $\forall c \in \Omega_c\cup \Omega_0, t \in \Omega_T, n \in \mathcal{G}, m \in \mathcal{D}, i\in \mathcal{B}$ and constraints (\ref{ramp_gen})-(\ref{shedding_relation}) hold $\forall c \in \Omega_c, t \in \Omega_T, m\in \mathcal{D}$.

The \textcolor{black}{annualized investment cost in VSR} is limited by the budget $A^{\max}$ in constraint (\ref{budget}). Constraints (\ref{p_bal})-(\ref{q_bal}) represent the active and reactive power balance at each bus. To guarantee that the power flow on the transmission line is zero when the line is in outage, a binary parameter $N_{kct}$ is introduced to denote the corresponding status of line $k$ for state $c$ under load level $t$ \cite{mybibb:investment_naps}. Note that we only consider transmission $N-1$ contingencies in this paper. However, the $N-k$ contingency can be easily implemented in the model with $N_{kct}$. Constraint (\ref{lS_normal}) ensures that no load shedding is allowed during base operating condition. The thermal limits of the transmission lines are considered in constraints (\ref{thermal_limit}). Constraints (\ref{Pg_limit1})-(\ref{Qg_limit1}) refer to the upper and lower bound for the active and reactive power production. The limits for the bus voltage magnitude and angle are enforced in (\ref{v_limit})-(\ref{theta_limit}). Constraints (\ref{vsr_limit}) represents the output range of the VSR. The reference bus angle is set to zero in constraint (\ref{ref_angle}).

Constraints (\ref{ramp_gen})-(\ref{shedding_relation}) are associated with the contingency states. In practical system, not all generators are capable of rescheduling after a certain contingencies. This constraint is enforced by (\ref{ramp_gen})-(\ref{g_fix}). Constraint (\ref{load_shedding}) ensures that the load shedding amount does not exceed the existing load. We assume that the power factor for the load is unchanged after the load curtailment, which is enforced by (\ref{shedding_relation}). 

\textcolor{black}{Note that the voltage magnitude $V$ and angle $\theta$ are state variables. However, in the OPF problem, these variables are usually treated as optimization variables \cite{mybibb:tcsc_wind}. Hence, the optimization variables of the planning model from (\ref{objective})-(\ref{shedding_relation}) are the elements in set $\splitatcommas{\Xi_{\text{OM}}=\{\Delta P^d_{mct},\Delta Q^d_{mct},\Delta P^{g,up}_{nct},\Delta P^{g,dn}_{nct},\delta_k,\theta_{ict},V_{ict},P^g_{nct},Q^g_{nct},x_{kct}^V\}}$.}

\section{Solution Approach}
\label{solution_approach}
The OPF is generally a non-convex and nonlinear problem which is hard to solve \cite{mybibb:scuc_dlr}. The introduction of the new variable ($\bm{x_k^V,\delta_k}$) to the optimization model makes the problem even more nonlinear. In addition, each variable in the planning model (\ref{objective})-(\ref{shedding_relation}) is usually associated with three dimensions, i.e., power system elements, states and load levels. The size of the planning model would dramatically increase with the system scale and considered operating states. Hence, the Generalized Benders Decomposition (GBD) \cite{mybibb:tcsc_ac,mybibb:gbd,mybibb:KTH_UC_uncertainty,mybibb:bd_scopf_m} is adopted to solve the proposed problem.

The complete optimization model is decomposed into a master problem and a number of subproblems. The master problem employs DC representation of the network to deal with the base operating condition for the three load levels. The subproblems exactly retain the AC characteristics of the network for all the considered operating states. The complicating variables between the master problem and subproblems are the active power generation $P^g_{n0t}$ and VSR installation $\delta_k$. The master problem and subproblems will be solved iteratively until the stopping criterion is satisfied. 

\begin{figure}[!htb]
	\centering
	\includegraphics[width=0.35\textwidth]{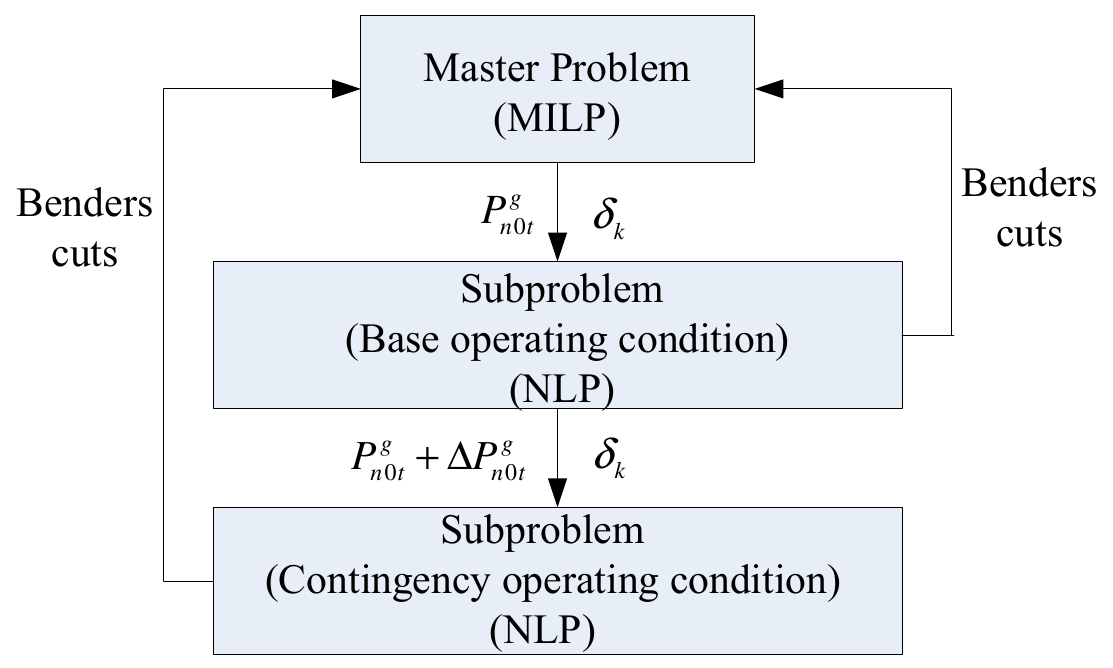}
	\caption{Flowchart of the proposed method.}
	\label{flowchart_ac}
\end{figure}

Note that the coupling constraint (\ref{ramp_gen}) between the base operating states and contingency operating states impedes the subproblem to be decomposed by states. To relieve the computational burden, a heuristic technique similar to the approach proposed in \cite{mybibb:KTH_UC_uncertainty} is leveraged. For each load level, the base operating state is solved at first in the subproblem and its active power generation will be the input for all the considered contingency states. Nevertheless, if necessary, a small number of contingencies will be solved with the base state in case of the severe contingencies. The flow chart of the proposed Benders algorithm is given in Fig. \ref{flowchart_ac}.

\subsection{Master Problem}
As mentioned at the beginning of this section, the master problem considers the base operating condition for the three load levels by using DC representation of the network. Note that the considered problem is still MINLP even for DC power flow (DCPF) model. The reformulation technique proposed in \cite{mybibb:tep_cvsr} is leveraged to transform the MINLP model into MILP model. For completeness, the reformulation technique is
illustrated in this section.
\subsubsection{Reformulation}
The power flow on the candidate transmission line $k$ to install VSR in DCPF can be expressed as:
\begin{align}
&P_{k}=(b_{k}+\delta_{k} b_{k}^{V})\theta_k,\ \ k\in \Omega_{V} \label{P_CVSR}  \\
&b_{k,V}^{\min} \le b_{k}^{V} \le b_{k,V}^{\max}, \ \ k\in \Omega_V   \label{bcvsr_range}
\end{align} 
In (\ref{P_CVSR}), $b_k^V$ is the susceptance change introduced by the VSR. It can be seen that the only nonlinearity lies in the trilinear term $\delta_{k}b_{k}^{V}\theta_{k}$ from (\ref{P_CVSR}). To eliminate the nonlinearity, a new variable $w_{k}$ is introduced:
\begin{equation}
w_{k}=\delta_{k}b_{k}^{V}\theta_{k}, \ \ k\in \Omega_{V} \label{wij}
\end{equation}
The constraint (\ref{P_CVSR}) is then modified to:
\begin{equation}
P_{k}=b_{k}\theta_k+w_{k}, \ \ k\in \Omega_{V} \label{PCVSR_wij}
\end{equation}
We multiply each side of the constraint (\ref{bcvsr_range}) with $\delta_{k}$ to yield:
\begin{equation}
\delta_{k} b_{k,V}^{\min} \le \frac{w_{k}}{\theta_{k}}=\delta_{k}b_{k}^{V} \le \delta_{k} b_{k,V}^{\max},\ \ k\in \Omega_{V}  \label{if_ineq}
\end{equation}
Depending on the sign of $\theta_{k}$, the inequality (\ref{if_ineq}) can be written as:
\begin{equation}
\left\lbrace 
\begin{aligned}
\delta_{k}\theta_{k} b_{k,V}^{\min} \le w_{k} \le \delta_{k}\theta_{k} b_{k,V}^{\max}, \ &\text{if} \ \theta_{k}>0  \\
w_{k}=0,\ \ \ \ \ \ \ \ \ \ \ \ \ \ \ \ \ \ \ \ \ \ \ \ \ &\text{if} \ \theta_{k}=0   \\
\delta_{k}\theta_{k} b_{k,V}^{\max} \le w_{k} \le \delta_{k}\theta_{k} b_{k,V}^{\min}, \ &\text{if} \ \theta_{k}<0 
\end{aligned} \right.
\end{equation}

The ``if" constraints can be formulated by introducing an additional binary variable $y_{k}$ and the big-M complementary constraints:
\begin{align}
&-M_{k}y_{k}+\delta_{k}\theta_{k}b_{k,V}^{\min} \le w_{k} \le \delta_{k}\theta_{k}b_{k,V}^{\max}+M_{k}y_{k},\ k\in \Omega_{V} \label{if_1} \\
&-M_{k}(1-y_{k})+\delta_{k}\theta_{k}b_{k,V}^{\max} \le w_{k}  \nonumber \\
&\ \ \ \ \ \ \ \ \ \ \ \ \ \ \ \  \le \delta_{k}\theta_{k}b_{k,V}^{\min}+M_{k}(1-y_{k}),\ k\in \Omega_{V} \label{if_2} 
\end{align}

Only one of these two constraints will be active during the optimization and the other one is a redundant constraint which is always satisfied because of the sufficiently large number $M_{k}$. However, the $M_{k}$ which is too large sometimes causes numerical problem. In this paper, $M_{k}$ is chosen to be equal to ($\max\{|b_{k,V}^{\min}|,|b_{k,V}^{\max}|\} \cdot \theta_{k}^{\max}$). 

In constraints (\ref{if_1}) and (\ref{if_2}), there still exists a bilinear term $\delta_{k}\theta_{k}$ which is the product of a binary variable and a continuous variable. We introduce another variable $z_{k}$ and use the standard linearization method to find: 
\begin{align}
&-\delta_{k}\theta_{k}^{\max} \le z_{k} \le \delta_{k}\theta_{k}^{\max},\ \ k \in \Omega_{V} \label{z1} \\
&\theta_{k}-(1-\delta_{k})\theta_{k}^{\max} \le z_{k} \le \theta_{k}+(1-\delta_{k})\theta_{k}^{\max},\ \ k\in \Omega_{V}  \label{z2}
\end{align}

Then the constraint (\ref{if_1}) and (\ref{if_2}) can be written as:
\begin{align}
&-M_{k}y_{k}+z_{k}b_{k,V}^{\min} \le w_{k} \le z_{k}b_{k,V}^{\max}+M_{k}y_{k},\ \ k\in \Omega_{V} \label{if3} \\
&-M_{k}(1-y_{k})+z_{k}b_{k,V}^{\max} \le w_{k}  \nonumber \\
&\ \ \ \ \ \ \ \ \ \ \ \ \ \ \ \ \ \ \ \ \  \le z_{k}b_{k,V}^{\min}+M_{k}(1-y_{k}),\ \ k\in \Omega_{V} \label{if4} 
\end{align}

The original MINLP model (\ref{P_CVSR})-(\ref{bcvsr_range}) has been transformed into an MILP model (\ref{PCVSR_wij}), (\ref{z1})-(\ref{if4}).   

\subsubsection{Master Problem Formulation}
With the reformulation, the master problem is formulated as follows:

\begin{align}
&\min_{\Xi_{\text{MP}}}\ Z_{down}^{(\nu)}=\sum_{t\in \Omega_T}\pi_{0t}\sum_{n \in \mathcal{G}}{a^g_nP_{n0t}^{g^{(\nu)}}}+\sum_{k\in \Omega_V}A^I_k\delta_k^{(\nu)} \nonumber \\
&\ \ \ \ \ \ \ \ \ \ \ \ \ \ \ \  +\sum_{c \in \Omega_c\cup \Omega_0 }\sum_{t \in \Omega_T}\alpha_{ct}^{(\nu)}   \label{master_obj} \\
& \text{subject to:}    \nonumber  \\
&(\ref{budget}), (\ref{Pg_limit1}),(\ref{theta_limit}),(\ref{ref_angle}), (\ref{PCVSR_wij}), (\ref{z1})-(\ref{if4}) \ \ \text{and}  \nonumber \\
&\sum_{n \in \mathcal{G}_i}P^{g^{(\nu)}}_{n0t}-\sum_{m \in \mathcal{D}_i}P^d_{m0t}=\sum_{k\in \Omega_{L}^{i}} P^{(\nu)}_{k0t}(\bm{\theta,w}) \label{p_bal_dc}   \\
& -S_{k0t}^{\max} \le P^{(\nu)}_{k0t}(\bm{\theta,w}) \le S_{k0t}^{\max}  \label{Smax} \\
&\alpha_{ct}^{(\nu)} \ge \alpha_{down}, \ \ c \in \Omega_c \cup \Omega_0  \label{accerlerate}  
\end{align}
\begin{align} 
&\alpha_{0t}^{(\nu)} \ge Z^{(l)}_{0t}+\sum_{n\in \mathcal{G}}\mu_{n0t}^{(l)}(P_{n0t}^{g^{(\nu)}}-P_{n0t}^{g^{(l)}}) \nonumber\\
&\ \ \ \  +\sum_{k\in \Omega_V}\beta_{k0t}^{(l)}(\delta_k^{(\nu)}-\delta_k^{(l)}), \ l=1,\cdots,\nu-1  \label{bender_cut_base}  \\
&\alpha_{ct}^{(\nu)} \ge Z^{(l)}_{ct}+\sum_{k\in \Omega_V}\beta_{kct}^{(l)}(\delta_k^{(\nu)}-\delta_k^{(l)}), \nonumber \\
&\ \ \ \ \ \ \ \ \ \ \ \ \ \ \ \ \ c \in \Omega_c, \ l=1,\cdots,\nu-1  \label{bender_cut_c} 
\end{align}
Constraints (\ref{master_obj})-(\ref{bender_cut_c}) hold $\forall t \in \Omega_T, n\in \mathcal{G}, i\in \mathcal{B}, k \in \Omega_L$.

All the variables are subject to Benders iteration parameter $\nu$. Constraint (\ref{p_bal_dc})-(\ref{Smax}) refer to the active power balance and thermal limit for the transmission line in the DC power flow respectively. A lower bound $\alpha_{down}$ is imposed on $\alpha_{ct}$ in (\ref{accerlerate}) to accelerate the convergence speed \cite{mybibb:KTH_UC_uncertainty}. Constraints (\ref{bender_cut_base})-(\ref{bender_cut_c}) are Benders cuts. Note that for each iteration, one cut is generated per operating state and load level. This is proved to be another trick to improve the convergence of the Benders algorithm \cite{mybibb:multicut}. The optimization variables of the master problem are those in the set $\Xi_{\text{MP}}=\{\theta_{i0t},P^{g}_{n0t},\delta_k,y_{k0t},z_{k0t},w_{k0t},\alpha_{ct}\}$.

\subsection{Subproblem}
With $\delta_k$ and $P^g_{n0t}$ from the master problem, each subproblem becomes a nonlinear and continuous program. The AC characteristics should be retained for both the base and contingency operating states. 

For the base operating conditions, i.e., $c \in \Omega_0$, the subproblem is formulated as:
\begin{align}
&\min_{\Xi_{\text{SP1}}} Z_{0t}^{(\nu)}=\pi_{0t}(\sum_{n \in \mathcal{G}}{a^g_n \Delta P_{n0t}^{g^{(\nu)}}}+\sum_{i\in \mathcal{B}}h_i^p(s_{i0t,1}^{p^{(\nu)}}+s_{i0t,2}^{p^{(\nu)}}) \nonumber \\
&\ \ \ \ \ \ \ +\sum_{i\in \mathcal{B}}h_i^q(s_{i0t,1}^{q^{(\nu)}}+s_{i0t,2}^{q^{(\nu)}}))  \label{sub_obj_base}  \\
&\text{subject to}   \nonumber \\
&(\ref{thermal_limit}),(\ref{Qg_limit1})-(\ref{ref_angle}) \ \ \text{and} \nonumber \\  
&\sum_{n \in \mathcal{G}_i}(P^{g^{(\nu)}}_{n0t}+\Delta P^{g^{(\nu)}}_{n0t})-\sum_{m \in \mathcal{D}_i}P^d_{m0t}  \nonumber \\
&\ \ \ \ \ \ \ \ +s^{p^{(\nu)}}_{i0t,1}-s^{p^{(\nu)}}_{i0t,2}=\sum_{k\in \Omega_{L}^i } P^{(\nu)}_{k0t}(\bm{\theta,V,x^V,\delta}) \label{p_bal_sp1}  \\
&\sum_{n \in \mathcal{G}_i}Q^{g^{(\nu)}}_{n0t}-\sum_{m \in \mathcal{D}_i}Q^d_{m0t}+s^{q^{(\nu)}}_{i0t,1}-s^{q^{(\nu)}}_{i0t,2} \nonumber \\
&\ \ \ \ \ \ \ \ \ \ \ \ \ \ \ \ \ \ \ \ \ \ \ \ \  =\sum_{k\in \Omega_{L}^i } Q^{(\nu)}_{k0t}(\bm{\theta,V,x^V,\delta}) \label{q_bal_sp1}  \\
&P^{g,\min}_{n0t}\le P^{g^{(\nu)}}_{n0t}+\Delta P^{g^{(\nu)}}_{n0t} \le P^{g,\max}_{n0t}  \label{pg_limit_base}  \\
&s_{i0t,1}^{p^{(\nu)}} \ge 0, \ s_{i0t,2}^{p^{(\nu)}} \ge 0,\ s_{i0t,1}^{q^{(\nu)}} \ge 0, \ s_{i0t,2}^{q^{(\nu)}} \ge 0    \label{slack_range}   \\
&P^{g^{(\nu)}}_{n0t}=\hat{P}^g_{n0t} \ \ \ \ : \mu_{n0t}^{(\nu)}   \label{fix_P} \\
&\delta_k^{(\nu)}=\hat{\delta}_k \ \ \ \ : \beta_{k0t}^{(\nu)}     \label{fix_de} 
\end{align} 
Constraints (\ref{sub_obj_base})-(\ref{fix_de}) hold $\forall t \in \Omega_T, n\in \mathcal{G}, i\in \mathcal{B}, k \in \Omega_L$.

The optimization variables of the base operating condition subproblem are those in the set $\splitatcommas{\Xi_{\text{SP1}}=\{\theta_{i0t},V_{i0t},P^g_{n0t},Q^g_{n0t},\Delta P^g_{n0t},s^p_{i0t,1},s^p_{i0t,2},s^q_{i0t,1},s^q_{i0t,2},\delta_k \}}$. In (\ref{sub_obj_base}), $\Delta P^g_{n0t}$ is the active power adjustment shifting from DCPF to ACPF. We introduce four slack variables $s_{i0t,1}^{p^{(\nu)}},s_{i0t,2}^{p^{(\nu)}},s_{i0t,1}^{q^{(\nu)}},s_{i0t,2}^{q^{(\nu)}}$ to ensure that the subproblems are always feasible, whose values are penalized in the objective function with sufficiently large constants $h^p_i$ and $h_i^q$. The objective function is to minimize the cost for the active power adjustment and the possible violations. Constraints (\ref{fix_P})-(\ref{fix_de}) fix the complicating variables to the values from the master problem. $\mu^{(\nu)}_{n0t}$ and $\beta^{(\nu)}_{n0t}$ are the sensitivities associated with these two constraints. 

According to the heuristic approach mentioned at the beginning of this section, the active power generation from the subproblem under base operating condition is the input for the subproblem under contingency operating condition. Hence, the subproblem for the contingency states can be formulated as:
\begin{align}
&\min_{\Xi_{\text{SP2}}} Z_{ct}^{(\nu)}=\pi_{ct}(C_{ct}+\sum_{i\in \mathcal{B}}h_i^p(s_{ict,1}^{p^{(\nu)}}+s_{ict,2}^{p^{(\nu)}}) \nonumber \\
&\ \ \ \ \ \ \ \ \ \ +\sum_{i\in \mathcal{B}}h_i^q(s_{ict,1}^{q^{(\nu)}}+s_{ict,2}^{q^{(\nu)}}))  \label{sub_obj_cont}  \\
&\text{subject to}   \nonumber \\
&(\ref{thermal_limit})-(\ref{ref_angle}),(\ref{ramp_up})-(\ref{ramp_dn}),(\ref{load_shedding})-(\ref{shedding_relation}) \ \ \text{and} \nonumber \\  
&\sum_{n \in \mathcal{G}_i}P^{g^{(\nu)}}_{nct}-\sum_{m \in \mathcal{D}_i}(P^d_{mct}-\Delta P^{d^{(\nu)}}_{mct})  \nonumber \\
&\ \ \ \ \ \ \ +s^{p^{(\nu)}}_{ict,1}-s^{p^{(\nu)}}_{ict,2}=\sum_{k\in \Omega_{L}^i } N_{kct}P^{(\nu)}_{kct}(\bm{\theta,V,x^V,\delta}) \label{p_bal_sp2}  \\
&\sum_{n \in \mathcal{G}_i}Q^{g^{(\nu)}}_{nct}-\sum_{m \in \mathcal{D}_i}(Q^d_{mct}-\Delta Q^{d^{(\nu)}}_{mct})  \nonumber \\
&\ \ \ \ \ \ \ +s^{q^{(\nu)}}_{ict,1}-s^{q^{(\nu)}}_{ict,2}=\sum_{k\in \Omega_{L}^i } N_{kct}Q^{(\nu)}_{kct}(\bm{\theta,V,x^V,\delta}) \label{q_bal_sp2}  \\
&P_{nct}^{g^{(\nu)}}=P_{n0t}^{g^{(\nu)}}+\Delta P^{g^{(\nu)}}_{n0t}+\Delta P_{nct}^{g,up^{(\nu)}}-\Delta P_{nct}^{g,dn^{(\nu)}}, \nonumber \\
&\ \ \ \ \ \ \ \ \ \ \ \ \ \ \ \ \ \ \ \ \ \ \ \ \ \ \ \ \ \ \ \ \ \ \ \ \ \ \ \ \ \ \ \ \ \ \ \ \ \ \ \ \ n\in \mathcal{G}_{re}    \label{ramp_gen_sub2} \\
&P_{nct}^{g^{(\nu)}}=P_{n0t}^{g^{(\nu)}}+\Delta P^{g^{(\nu)}}_{n0t}, n\in \mathcal{G}\backslash\mathcal{G}_{re}  \label{g_fix_sub2}  \\
&s_{ict,1}^{p^{(\nu)}} \ge 0, \ s_{ict,2}^{p^{(\nu)}} \ge 0,\ s_{ict,1}^{q^{(\nu)}} \ge 0, \ s_{ict,2}^{q^{(\nu)}} \ge 0    \label{slack_range_sub2}   \\
&\delta_k^{(\nu)}=\hat{\delta}_k \ \ \ \ : \beta_{kct}^{(\nu)}     \label{fix_de_sub2} 
\end{align} 
Constraints (\ref{sub_obj_cont})-(\ref{fix_de_sub2}) hold $\forall c \in \Omega_c,t \in \Omega_T, n\in \mathcal{G}, i\in \mathcal{B}, k \in \Omega_L,m \in \mathcal{D}$.

The optimization variables of the subproblem under contingency states are those in the set $\splitatcommas{\Xi_{\text{SP2}}=\{\theta_{ict},V_{ict},P^g_{nct},Q^g_{nct},\Delta P^{g,up}_{nct},\Delta P^{g,dn}_{nct},\Delta P^d_{mct},\Delta Q^d_{mct},\delta_k,s^p_{ict,1},s^p_{ict,2},s^q_{ict,1},s^q_{ict,2}\}}$. The objective function is to minimize the operating cost associated with the contingency states and the possible violations. The sensitivities used to construct Benders cut are generated by (\ref{fix_de_sub2}).

At each iteration $\nu$, the upper bound of the objective function for the original problem can be calculated as follows:
\begin{align}
&Z_{up}^{(\nu)}=\sum_{t \in \Omega_T}Z_{0t}^{(\nu)}+\sum_{t \in \Omega_T}\sum_{c \in \Omega_c}Z_{ct}^{(\nu)} \nonumber \\
&\ \ \ \ \ \ \ \ \ \ +\sum_{t\in \Omega_T}\pi_{0t}\sum_{n \in \mathcal{G}}{a^g_n\hat{P}_{n0t}^{g}}+\sum_{k\in \Omega_V}A^I_k\hat{\delta}_k
\end{align} 

The iteration procedure will end until all the slack variables are zero and the difference between the upper bound and lower bound for the objective function is within a predefined tolerance $\epsilon$:
\begin{equation}
\frac{|Z_{up}^{(\nu)}-Z_{down}^{(\nu)}|}{|Z_{up}^{(\nu)}|}\le \epsilon
\end{equation} 

\begin{figure*}[!htb]
	\centering
	\includegraphics[width=0.85\textwidth]{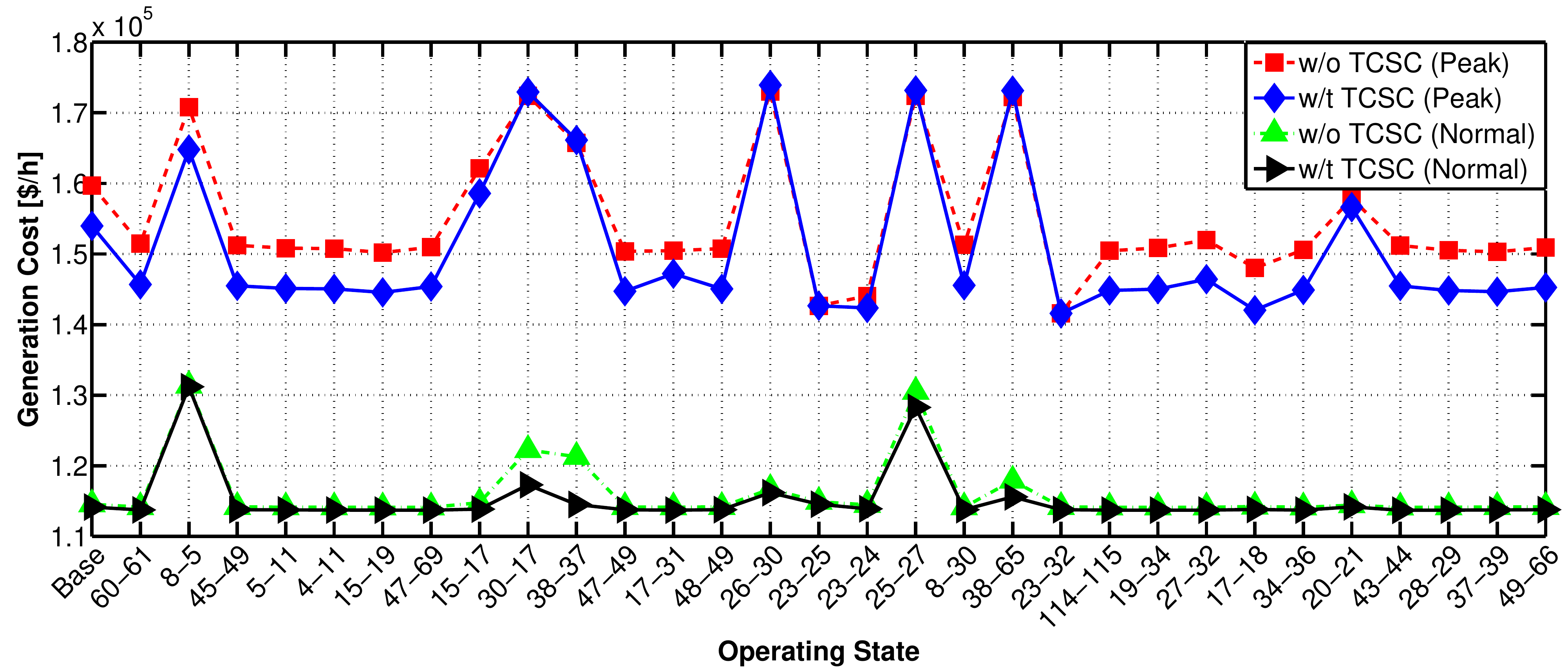}
	\caption{Hourly generation cost for peak and normal load level of IEEE 118-bus system.}
	\label{oper_cost_118}
\end{figure*}

\section{Numerical Case Studies}
\label{case_study}
We test our proposed planning model and solution approach on the IEEE 118-bus system. The data for the system is from the MATPOWER software \cite{mybibb:MATPOW}. Since only one snapshot of load is provided for this test system in \cite{mybibb:MATPOW}, we treat that load as the normal load level. The peak load level is 20\% higher than the normal load level and the low load level is 20\% less than the normal load level. The computer used for all the simulations has an Inter Core(TM) i5-2400M CPU @ 2.30 GHz with 4.00 GB of RAM. The complete model is programmed in GAMS \cite{mybibb:GAMS}. The MILP master problem is solved by CPLEX \cite{mybibb:CPLEX} and the NLP subproblems are solved by IPOPT \cite{mybibb:ipopt}. 

We consider the allocation strategy for the TCSC in this paper. The compensation rate for the TCSC is allowed to vary from -70\% to 20\% of its corresponding reactance \cite{mybibb:TCSC_cost_recovery}. \textcolor{black}{According to \cite{mybibb:GA_FACT_market}, the annual investment cost of TCSC is converted from its total investment cost with the yearly interest rate and life span using equation (\ref{cost_convert}). In this study, the interest rate $d$ is assumed to be 5\% and the life span $LT$ is selected to be 5 years \cite{mybibb:TCSC_cost_recovery}.} 
\begin{equation}
\textcolor{black}{\tilde{A}^I_k={A}^I_k \cdot \frac{d(1+d)^{LT}}{(1+d)^{LT}-1}}  \label{cost_convert}
\end{equation}    
The budget $A^{\max}$ is assumed to be 3 M\$. The algorithm tolerance $\epsilon$ is set to be 0.2\%. The constant $\alpha_{down}$ in (\ref{accerlerate}) is selected to be $-10^{10}$ through a trial-and-error process. \textcolor{black}{The duration of peak, normal and low load level are assumed to be 2190 h, 4380 h and 2190 h, respectively. In addition, based on the typical values of line outage rate provided in \cite{mybibb:force_outage_1979}, we assume that the forced outage rate for the contingency branch is 0.1\%. Table \ref{duration} gives the number of operating hours in each state.} Finally, the planning criterion for the base and contingency operating condition is provided in Table \ref{criterion}. 
\begin{table}[!htb]
	\centering
	\caption{\textcolor{black}{Duration of Each Operating State}}
	\label{duration}
	\begin{tabular}{c c c c}
		\hline
		&Peak&Normal&Low  \\
		\hline
		\# of hours for base state&2124.3&4248.6&2124.3   \\
		\hline
		\# of hours for each contingency state&2.19&4.38&2.19   \\
		\hline   
	\end{tabular}
\end{table}

\begin{table}[!htb]
	\centering
	\caption{Planning Criterion for the Test System}
	\label{criterion}
	\begin{tabular}{c c c }
		\hline
		&Base&$N-1$ Contingency  \\
		\hline
		Voltage (p.u.)&$0.94 \le V \le 1.06$& $0.9 \le V \le 1.1$  \\
		\hline 
		Thermal Limits &$S^{\max}_k$&$1.1S^{\max}_k$   \\
		\hline
	\end{tabular}
\end{table}

\subsection{IEEE 118-Bus System}
The IEEE 118-bus system has 118 buses, 177 transmission lines, 9 transformers and 19 generators. The total active and reactive load for the peak load level are 4930 MW and 1661 MVar. The active and reactive power generation capacity are 6466 MW and 6325 MVar. The thermal limits for the transmission lines are decreased artificially to create congestions. We consider 30 contingencies based on the congestion severity \cite{mybibb:automatic_contingency} so the number of operating states is 93. Moreover, we leverage the sensitivity approach \cite{mybibb:PSO_SQP_FACTS} to select 30 candidate locations to install TCSC. \textcolor{black}{The complete procedures are given below:}
\begin{enumerate}
	\item [1)] 	\textcolor{black}{Run an OPF for each operating state without TCSC using IPOPT as the solver. Note that the coupling constraints between the base state and contingency state are not considered. In addition, each line reactance is treated as an optimization variable and constraint (\ref{fix_x}) is included in the OPF model, i.e., fix the line reactance to its original value:}
	\begin{equation}
	\tilde{x}_{kct}=x_k  \label{fix_x}	
	\end{equation}
     \item [2)]  \textcolor{black}{Obtain the sensitivity ($\lambda_{kct}$) of the operation cost with respect to the change of line reactance in each state, i.e., value of the dual variable associated with constraint (\ref{fix_x}).}
     \item [3)] \textcolor{black}{Compute the weighted sensitivity ($\bar{\lambda}_k$) of branch $k$ by equation (\ref{weighted_sens}):}
     \begin{equation}
     \bar{\lambda}_k=\sum_{t\in \Omega_T} \sum_{c \in \Omega_c}\pi_{ct}|\lambda_{kct}x_k|  \label{weighted_sens}
     \end{equation}
     \item [4)] \textcolor{black}{Sort $\bar{\lambda}_k$ in a descending order and select the first 30 lines as candidate locations for TCSC.}
 \end{enumerate}

\begin{table}[!htb]
	\centering
	\caption{Annual Planning Cost Comparison for IEEE 118-Bus System}
	\label{annual_saving_118}
	\begin{tabular}{c c c}
		\hline
		\multirow{2}{*}{Cost Category}&\multicolumn{2}{c}{Annual Cost [million \$]}     \\
		\cline{2-3}
		&w/o TCSC&w/t TCSC      \\
		\hline
		Generation cost in normal state&1013.85&999.97   \\
		\hline
		Generation cost in contingency&31.23&30.85     \\
		\hline
		Rescheduling cost&0.84&0.76           \\
		\hline
		Load shedding cost&17.81&13.46      \\
		\hline
		Investment on TCSC&-&2.84     \\
		\hline
		Total cost&1063.72&1047.88      \\
		\hline
		&&15-33, 17-18  \\
		&&28-29, 24-72   \\
		TCSC locations ($i-j$)&-&20-21, 74-75  \\
		&&40-42, 22-23  \\
		&&35-37, 37-39  \\
		\hline
		Computational Time [s]&44.47 &437.83  \\
		\hline
	\end{tabular}
\end{table}

The planning model suggests that 10 transmission lines will be installed with TCSCs. Table \ref{annual_saving_118} compares the annual planning cost for the case without and with TCSCs. The second row provides the generation cost in the base operating states, i.e., the second term in (\ref{objective}). The generation cost in contingency states, i.e., the first term in (\ref{obj_cont}), is given in the third row. The rescheduling and load shedding cost in the contingency states are provided in the fourth and fifth row. The sixth row shows the annual investment cost on TCSC. The annual total planning cost is provided in the seventh row. The eighth row gives the selected TCSC locations. The computational time for the planning models is provided in the last row. It can be seen that the installation of the TCSCs decreases the cost in all the categories. Although the investment on TCSCs costs 2.84 M\$, the annual saving is about 15.84 M\$. The computational time for the case considering TCSC is about 437.83 s for the IEEE-118 bus system.

\begin{figure}[!htb]
	\centering
	\includegraphics[width=0.45\textwidth]{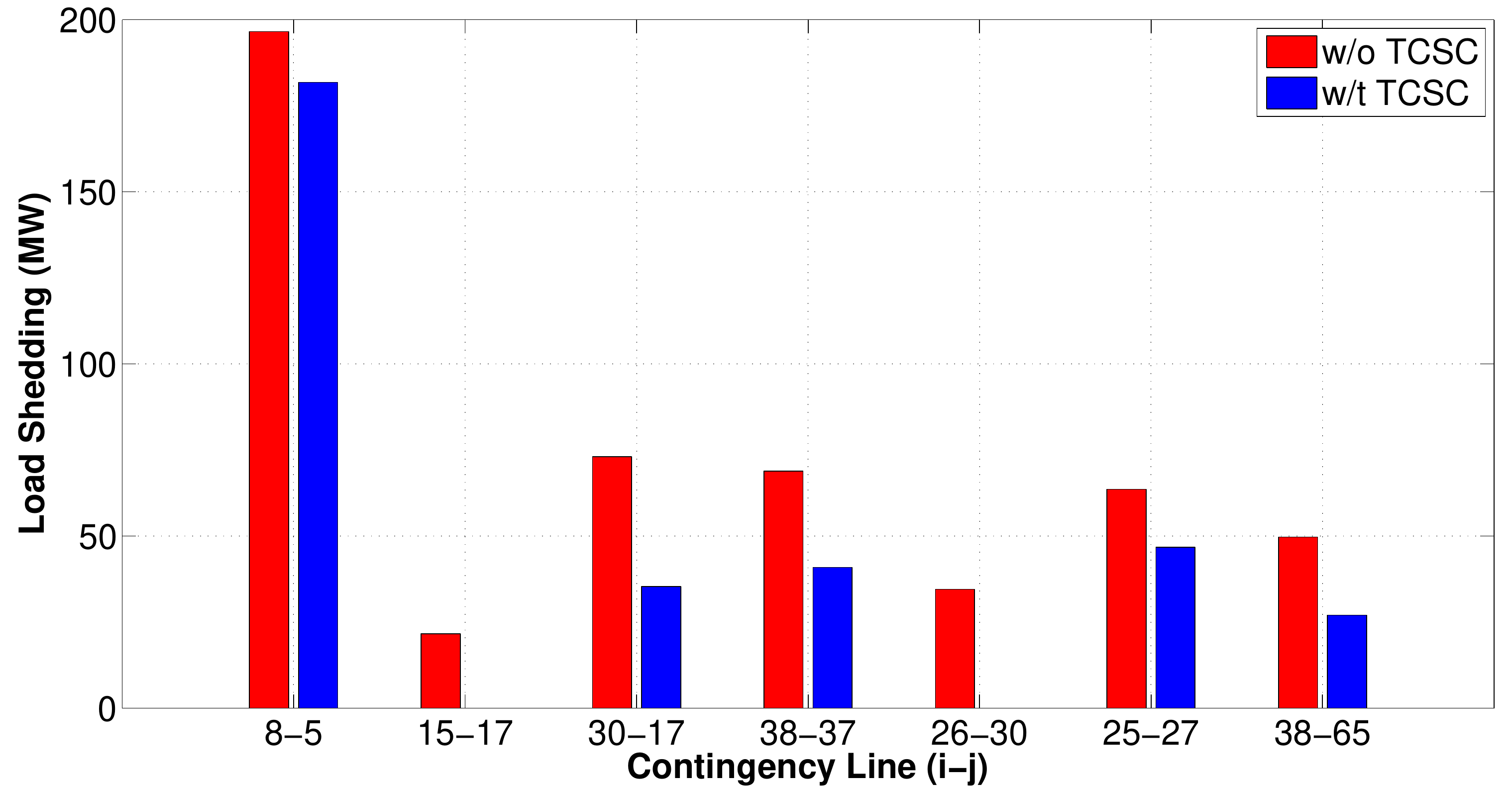}
	\caption{Load shedding amount under different contingencies for peak load level.}
	\label{LS}
\end{figure}

Fig. \ref{oper_cost_118} shows the hourly generation cost for each states under the peak and normal load level and Fig. \ref{LS} provides the load shedding amount under the peak load level for the 7 contingencies in which the load shedding exists. From Fig. \ref{oper_cost_118}, it can be seen that the generation cost is reduced in the majority of operating states under both the peak and normal load level. For instance, the generation cost for contingency (60-61) is 151447 \$/h under the peak load level without TCSC. The cost is decreased by 5761 \$/h with TCSCs. Under the peak load level, the generation cost for contingency (26-30), (25-27) and (38-65) with TCSCs are slightly higher than that without TCSC. Nevertheless, significant load curtailment reductions can be observed for these three contingencies from Fig. \ref{LS}. Therefore, the total operating cost for the three contingencies are still cheaper by the installation of TCSCs. Under the normal load level, the installation of TCSCs decreases the generation cost for all the states. However, the cost reductions for most states are not as much as that under peak load level. The cost reductions are mainly due to the congestion relief so that more load can be covered by cheap generators. From Fig. \ref{LS}, it can be seen that the load shedding reductions occur in all the 7 contingencies by installing TCSCs. The load shedding for contingency (15-17) and (26-30) are completely eliminated. The largest load shedding reduction is for contingency (30-17), the amount of the load curtailment decreases from 73.05 MW to 35.35 MW. 

\begin{figure}[!htb]
	\centering
	\includegraphics[width=0.4\textwidth]{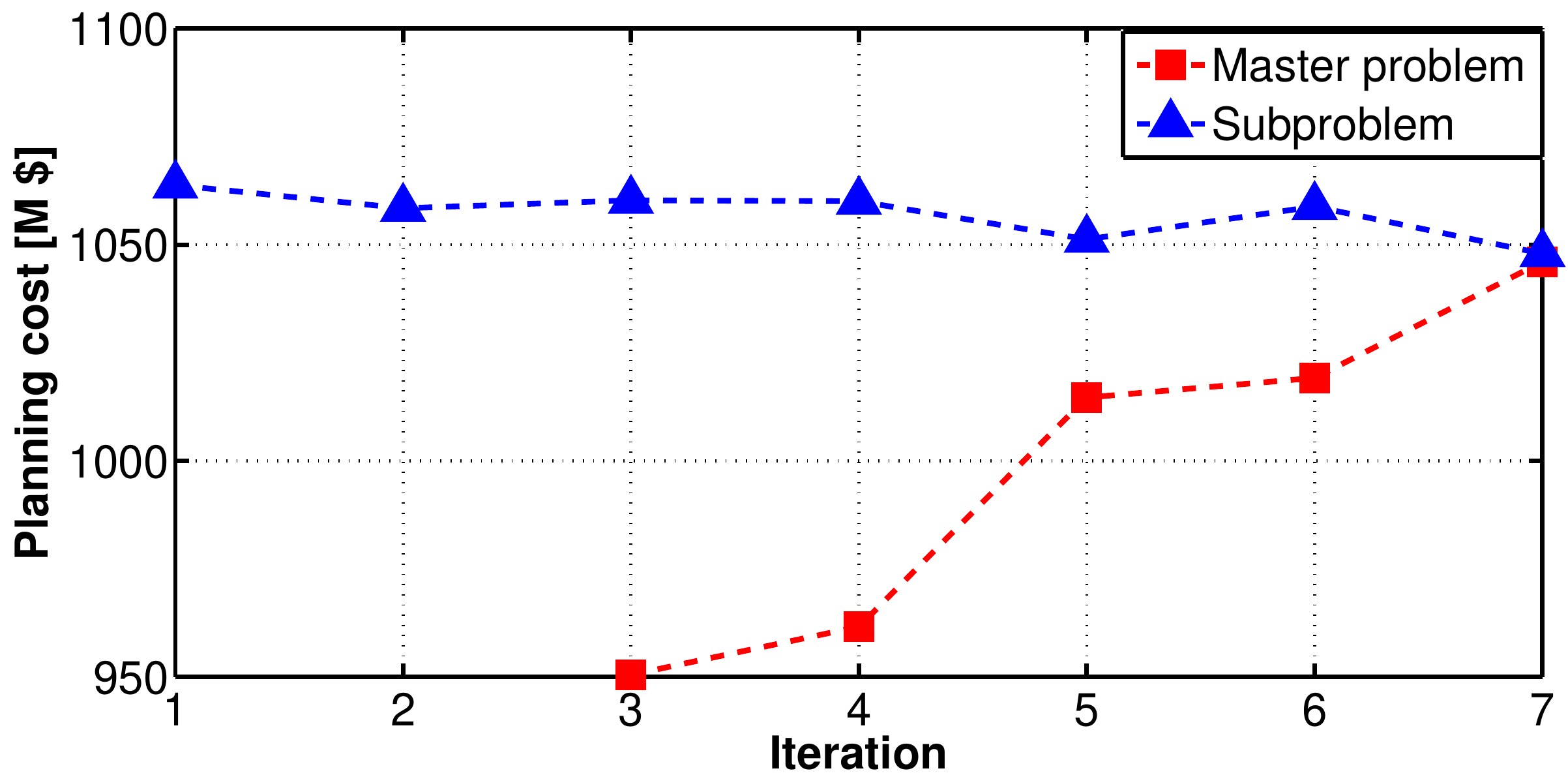}
	\caption{Evolution of the proposed Benders algorithm for IEEE 118-bus system.}
	\label{iteration}
\end{figure}

Fig. \ref{iteration} illustrates the iteration process for the proposed algorithm. Note that the objective values of lower bound for the first two iterations are negative so they are not plotted to make the figure more readable. It can be seen that the difference between the upper bound $Z_{up}$ and the lower bound $Z_{down}$ is within the tolerance at iteration 7, where the algorithm converges.  

\subsection{\textcolor{black}{Computational Issues}}
\textcolor{black}{As mentioned in the introduction, the proposed planning model is a large scale MINLP problem which is beyond the capability of the existing commercial solvers. For comparison, we leverage BONMIN \cite{mybibb:bonmin} to solve the complete model, i.e., (\ref{objective})-(\ref{shedding_relation}). The time limit is set to be six hours. The model size reported by BONMIN is provided in Table \ref{model_size}. After six hours run on the personal computer, the message returned by BONMIN is \textsl{no feasible solution found.} The results further demonstrate the advantages of our proposed decomposition algorithm for the TCSC allocation problem in large scale power networks.}
\begin{table}[!htb]
	\centering
	\caption{\textcolor{black}{Model Size of the Proposed Planning Model}}
	\label{model_size}
	\begin{tabular}{c c c c}
		\hline
		\# of binary&\# of continuous&\# of equality&\# of inequality   \\
		variables&variables&constraints&constraints  \\
		\hline
		30&124908&107787&34597 \\
		\hline

	\end{tabular}
\end{table}

\section{Conclusion and Discussion}
\label{conclusions}
In this paper, a planning model to optimally allocate VSRs in the transmission network considering AC constraints and a series of transmission $N-1$ contingencies is proposed. Originally, the planning model is a large scale MINLP model that is generally intractable. To relieve the computational burden, Generalized Benders Decomposition is utilized to separate the model into master problem and a number of subproblems. The numerical results based on the IEEE-118 bus system validate the performance of the proposed approach. In addition, the simulation results show that the appropriate installation of VSRs decreases the operating cost in both the base and contingency conditions and allows reduced planning cost.

\textcolor{black}{Note that in practice, there can be many factors to be considered rather than planning cost and $N-1$ margins, e.g., policy issues, location constraints, space restrictions, land cost, etc. In such cases, rather than providing the optimal solution, it is desired to provide a set of optimal solutions (with rankings). It is up to the decision makers to synthesize all impact factors and select the best solution that works for them, which typically involves compromise. Our future work is to leverage the multi-objective optimization framework to determine the locations of VSR that co-optimize the planning cost and the transmission margins.}


%




\ifCLASSOPTIONcaptionsoff
  \newpage
\fi



%
\bibliographystyle{IEEEtran}
\bibliography{IEEEabrv,mybibb}

%

\begin{IEEEbiographynophoto}{Xiaohu Zhang}
	(S'12) received the B.S. degree in electrical engineering from Huazhong University of Science and Technology, Wuhan, China, in 2009, the M.S. degree in electrical engineering from Royal Institute of Technology, Stockholm, Sweden, in 2011, and the Ph.D. degree in electrical engineering at The University of Tennessee, Knoxville, in 2017. 
	
	Currently, he works as a power system engineer at GEIRI North America, San Jose, CA, USA. His research interests are power system operation, planning and stability analysis.
\end{IEEEbiographynophoto}

\begin{IEEEbiographynophoto}{Chunlei Xu}
	received the B.S. degree in electrical
	engineering from Shanghai Jiao Tong University,
	Shanghai, China, in 1999. He currently leads the
	Dispatching Automation department at Jiangsu
	Electrical Power Company in China. His research
	interests include power system operation and control
	and WAMS.
\end{IEEEbiographynophoto}

\begin{IEEEbiographynophoto}{Di Shi}
	(M’12, SM’17) received the Ph.D. degree in
	electrical engineering from Arizona State University,
	Tempe, AZ, USA, in 2012. He currently leads the
	Advanced Power System Analytics Group at GEIRI
	North America, San Jose, CA, USA. He has
	published over 50 journal and conference papers and
	hold 13 US patents/patent applications.
\end{IEEEbiographynophoto}

\begin{IEEEbiographynophoto}{Zhiwei Wang}
	received the B.S. and M.S. degrees in
	electrical engineering from Southeast University,
	Nanjing, China, in 1988 and 1991, respectively. He is
	President of GEIRI North America, San Jose, CA,
	USA. His research interests include power system
	operation and control, relay protection, power system
	planning, and WAMS.
\end{IEEEbiographynophoto}

\begin{IEEEbiographynophoto}{Qibing Zhang}
	received the B.S. degree in electrical
	engineering from Zhejiang University, Zhejiang, China,
	in 2007, and the M.S. degree in electrical engineering
	from Shanghai Jiao Tong University, Shanghai, China,
	in 2010. His research interests include power system
	operation and control and relay protection.
\end{IEEEbiographynophoto}

\begin{IEEEbiographynophoto}{Guodong Liu}
	(S’10, M’14) received his B.S. in electrical engineering from Shandong University, Jinan, China, in 2007, M.S. in electrical engineering from Huazhong University of Science and Technology, Wuhan, China, in 2009, and Ph.D. in electrical engineering from the University of Tennessee, Knoxville in 2014. He is currently a R\&D staff in the Electrical and Electronic System Research Division at Oak Ridge National Laboratory where he leads projects on microgrid planning and operation, renewable energy integration and active distribution network management. His research interest includes power system economic operation and reliability, distribution management system, renewable energy and microgrids.
\end{IEEEbiographynophoto}

\begin{IEEEbiographynophoto}{Kevin Tomsovic}
	(F'07) received the BS from Michigan Tech. University, Houghton, in 1982, and the MS and Ph.D. degrees from University of Washington, Seattle, in 1984 and 1987, respectively, all in Electrical Engineering. He is currently the CTI Professor with the Department of Electrical Engineering and Computer Science, University of Tennessee, Knoxville, TN, USA, where he directs the NSF/DOE ERC, Center for Ultra-Wide-Area Resilient Electric Energy Transmission Networks (CURENT),  and previously served as the Electrical Engineering and Computer  Science Department  Head  from 2008  to  2013.  He was on  the  faculty of Washington State University, Pullman, WA, USA, from 1992 to 2008. He held the Advanced Technology for Electrical Energy Chair at National Kumamoto University, Kumamoto, Japan, from  1999 to 2000, and  was the NSF Program  Director with  the Electrical and Communications Systems Division of the Engineering Directorate from 2004 to 2006. He also held positions at National Cheng Kung University and National Sun Yat Sen University in Taiwan from 1988-1991 and the Royal Instituted of Technology in Sweden from 1991-1992. He is a Fellow of the IEEE.
\end{IEEEbiographynophoto}

\begin{IEEEbiographynophoto}{Aleksandar Dimitrovski}
	(SM'09) is an Associate Professor at the University of Central Florida, Orlando. Before joining UCF, he was a chief technical scientist at the Oak Ridge National
	Laboratory and a joint faculty at the University of Tennessee, Knoxville. In the past, he had been with Schweitzer Engineering Laboratories, and Washington State University, Pullman. He received his B.Sc. and Ph.D. degrees in electrical engineering with emphasis in power, and M.Sc. degree in applied computer sciences in Europe. His area of interest has been focused on modeling, analysis, protection, and control of uncertain power systems and, recently, on hybrid magnetic-electronic power control devices.
\end{IEEEbiographynophoto}






\enlargethispage{-3in}

\end{document}